\documentclass{article}

\usepackage{arxiv}
\usepackage{amsmath}
\usepackage[utf8]{inputenc} 
\usepackage[T1]{fontenc}    
\usepackage{booktabs}       
\usepackage{amsfonts}       
\usepackage{nicefrac}       
\usepackage{microtype}      
\usepackage[square,sort,comma,numbers]{natbib}
\usepackage{doi}
\usepackage{url}

\title{Remark on Fermat's Last Theorem}

\date{August 22, 2021}

\author{{Deriyenko I.I.} \\
	\texttt{ivan.deriyenko@gmail.com} \\
}

\hypersetup{
pdftitle={Remark on Fermat's Last Theorem},
pdfsubject={Remark on Fermat's Last Theorem},
pdfauthor={Deriyenko I.I.},
pdfkeywords={},
}

\begin{document}
\maketitle

\begin{abstract}
In this short article we do not prove Fermat's last theorem. We show that the number 2 is an exceptional number in this theorem.
\end{abstract}

Consider the well-known equation

\begin{equation}\label{eq1}
    X^n+Y^n=Z^n
\end{equation}

where $X,Y,Z,n\in\mathbb{N}$ and $n\neq1$. Let $X<Y$. Let's the following elementary transformations: 
\begin{equation*}
    \left(\frac{X}{Z}\right)^n+\left(\frac{Y}{Z}\right)^n=1
\end{equation*}

We take the following notation: $\frac{X}{Z}=x,\frac{Y}{Z}=y$, where $x,y\in\mathbb{Q}$. It's obvious that $0<x<y<1$, then we have:
\begin{equation*}
    y=\left(1-x^n\right)^{\frac{1}{n}}
\end{equation*}

take the indefinite integral of both parts:
\begin{equation*}
    \int y\left(x\right)dx=\int \left(1-x^n\right)^{\frac{1}{n}}dx
\end{equation*}

And we got integration of a differential binomial, in accordance with P.L.Chebyshev
\begin{equation}\label{eq2}
    I=\int x^a \left(\alpha+\beta x^b\right)^c dx
\end{equation}
is expressible by of the elementary function only in these three cases \cite{integration}:
\begin{equation*}
    \frac{a+1}{b}+c\in\mathbb{Z}
\end{equation*}
\begin{equation*}
    \frac{a+1}{b}\in\mathbb{Z}
\end{equation*}
\begin{equation*}
    c\in\mathbb{Z}
\end{equation*}
where $\mathbb{Z}$ is a whole number.

In our case we have:
\begin{equation*}
    a=0, b=n, c=\frac{1}{n}, \alpha=1, \beta=-1
\end{equation*}
so 
\begin{equation*}
    c=\frac{1}{n}\notin\mathbb{Z}
\end{equation*}
\begin{equation*}
    \frac{a+1}{b}=\frac{1}{n}\notin\mathbb{Z}
\end{equation*}
for $n\in\mathbb{N}, n\neq 1$
and
\begin{equation*}
    \frac{a+1}{b}+c=\frac{1}{n}+\frac{1}{n}=\frac{2}{n}\in\mathbb{Z}
\end{equation*}
at $n=2$ only.

If $n=2$ we have elementary function:
\begin{equation*}
    \int \sqrt{1-x^2}dx=\frac{1}{2}\left(x\sqrt{1-x^2}+\arcsin{x}\right)+C
\end{equation*}

In the other case $\left(n\neq2\right)$ we get from \ref{eq2} nonelementary functions \cite{integration}
\begin{equation*}
    I=\frac{1}{1+a}\alpha^{\frac{a+1}{b}}\beta^{-\frac{a+1}{b}}y^{\frac{1+a}{b}}F\left(\frac{a+1}{b};2-c;\frac{1+a+b}{b};y\right);y=\beta \frac{x^{b}}{\alpha}
\end{equation*}

in our case we have
\begin{gather*}
    \frac{1}{1+a}=1; \alpha^{\frac{a+1}{b}+c}=1^{\frac{2}{n}}=1;\beta^{-\frac{a+1}{b}}=(-1)^{-\frac{1}{n}};y^{\frac{1+a}{b}}=\left(-x^n\right)^{\frac{1}{n}}=-x \\
    \beta^{-\frac{a+1}{b}}\cdot y^{\frac{a+1}{b}}=\left(\beta^{-1}y\right)^{\frac{a+1}{b}}\left(-1^{-1}\left(-x^n\right)\right)^{\frac{1}{n}}=x \\
    I=xF\left(\frac{1}{n};2-\frac{1}{n};1+\frac{1}{n};-x^n\right)
\end{gather*}
where $F\left(p;q;r;z\right)=1+\frac{pq}{r}\cdot z+\frac{p\left(p+1\right)q\left(q+1\right)}{r\left(r+1\right)2!}\cdot z^2+\cdots$ --- hypergeometric function defined by Gauss series \cite{DLMF}. If $q=r$, then $2-\frac{1}{n}=1+\frac{1}{n}$, hence $n=2$ and $F$ transforms into a regular series (elementary function).

This shows the exclusivity of the number 2 in the equality \ref{eq1}, if $n=2$ we get an elementary function and in all the other cases we get a non-elementary function.

\bibliographystyle{unsrtnat}


\begin{thebibliography}{1}
    \bibitem{integration}
    \newblock Integration of differential binomial, 2014 \\
    \url{https://planetmath.org/IntegrationOfDifferentialBinomial}
    \bibitem{DLMF}
    \newblock DLMF 15.1 Special Notation, 2014 \\
    \url{https://dlmf.nist.gov/15.2#i}
\end{thebibliography}





\end{document}